\documentclass{amsart}

\usepackage{amssymb}

\theoremstyle{plain}

\numberwithin{equation}{section}

\newtheorem{theo}[equation]{Theorem}

\newtheorem*{teo}{Theorem A}
\newtheorem*{teob}{Theorem B}

 \newtheorem*{conj}{Conjecture C}

\newtheorem{question}[equation]{Question}

\newtheorem{lem}[equation]{Lemma}

\newcommand{\dl}{\operatorname{dl}}

\newcommand{\Ker}{\operatorname{Ker}}

\theoremstyle{definition}

\newtheorem{definition}[equation]{Definition}
\newtheorem{example}[equation]{Example}

\begin{document}	
\title{On conjugacy classes and derived length}

\author{Edith Adan-Bante}

\address{Department of Mathematical Science, Northern Illinois University, 
 Watson Hall 320
DeKalb, IL 60115-2888, USA}

\email{EdithAdan@illinoisalumni.org}


\subjclass{20d15}

\date{2009}
\begin{abstract} Let $G$ be a finite group and $A$, $B$ and $D$  be conjugacy classes of
$\,G$ with $D\subseteq AB=\{xy\mid x\in A, y\in B\}$. Denote by $\eta(AB)$ the number
of distinct conjugacy classes such that $AB$ is the union of those. Set 
${\bf C}_G(A)=\{g\in G\mid x^g=x \mbox{ for all } x\in A\}$.
If $AB=D$ then ${\bf C}_G(D)/({\bf C}_G(A)\cap{\bf C}_G(B))$ is an 
abelian group. If, in addition, $G$ is supersolvable, then the derived length of 
${\bf C}_G(D)/({\bf C}_G(A)\cap{\bf C}_G(B))$ is bounded above by $2\eta(AB)$. 
\end{abstract}

\maketitle

\begin{section}{Introduction}

 Let $G$ be a finite group, $A\subset G$ be a conjugacy class of G, i.e.
   $A=a^G=\{g^{-1} a g\mid g\in G\}$  for some $a$ in $G$.
Let $X$ be a normal subset of $G$, i.e. 
$X^g=\{ g^{-1} xg\mid x\in X\}=X$ for all $g\in G$.
  We can check that $X$ is  a union of 
  $n$ distinct conjugacy classes of $G$, for some integer $n>0$. Set
 $\eta(X)=n$.
 
  We can check that given any two conjugacy classes
 $A$ and $B$ of $G$, the product $A B=\{xy \mid x\in A  , y\in B\}$ of
 $A$ and $B$ is a normal subset of $G$. Then $\eta(AB)$ is
   the number of distinct conjugacy 
 classes of $G$ such that $AB$ is the union of those classes.
 
Denote by ${\bf C}_G(X)= \{g\in G \mid x^g=x \mbox{ for all }x\in X \}$ the centralizer of $X$ in $G$.  If $G$ is a
solvable group, denote by $\dl(G)$ the derived length of $G$. Let 
${\bf Z}(G)$ be the center of the group $G$.

In this note, we are exploring the relations between the structure of the group $G$ and 
the product $AB$ of  some conjugacy 
 classes $A$ and $B$ of $\,G$. More specifically,
 we are exploring the relation  between the derived length of some section of $G$ and
 properties of $AB$.

Given a finite solvable group $G$ and conjugacy classes $A$ and $B$ of $G$,
is there any relationship between the derived length of $G$ and 
 $\eta(AB)$? In general, the answer seems to be no. For instance
$A\{e\}=A$ for any finite group $G$ and any conjugacy class $A$ of $G$. 
 Thus $\eta(AB)$ may not give us information
 about $\dl(G)$, but it does give us a linear bound on the derived length of a 
 section of $G$, namely on the 
 section ${\bf C}_G(D)/( {\bf C}_G(A)\cap{\bf C}_G(B))$.

 \begin{teo} Let $G$ be a finite group and $A$, $B$ be normal subsets of $G$. 
 Then
 ${\bf C}_G(AB)/( {\bf C}_G(A)\cap{\bf C}_G(B))$ is abelian.
 \end{teo}
 
 Given any integer $m>0$,  we show in Example \ref{bunchextraspecials} 
 that there exists a nilpotent group $G$
 with conjugacy classes $A$, $B$ and $D$ such that
  $|{\bf C}_G(D)/( {\bf C}_G(A)\cap{\bf C}_G(B))|=m$ and $AB=D$. It follows that
  although 
  ${\bf C}_G(D)/( {\bf C}_G(A)\cap{\bf C}_G(B))$ is abelian, its order is unbounded and 
  even the number of distinct prime divisors is unbounded.

 \begin{teob}
 Let $G$ be a supersolvable group, $A$, $B$ and $D$ conjugacy classes of $\,G$ such that 
$D\subseteq AB$. Then 
$\dl ({\bf C}_G(D)/({\bf C}_G(A)\cap{\bf C}_G(B)))\leq 2\eta(AB)$.
\end{teob}

 We now mention the ``dual" situation for characters, where the ``dual" of a conjugacy
 class $A$ is an irreducible character $\chi$ and the ``dual" of the the Kernel 
 $\Ker(\chi)$ of $\chi$ is
 the centralizer ${\bf C}_G(A)$ of $A$.
 Let  $\Xi$ and $\Psi$
 be 
complex characters of $G$. Since any product of characters is a character,
$\Xi\Psi$ is a character of $G$. Thus it can be written as
an integral linear combination of irreducible characters of $G$.
 Let $\eta(\Xi\Psi)$ be the number of distinct
irreducible constituents
of the character $\Xi\Psi$. In Theorem A of \cite{derivedcharacters2}, it is proved that
there exist universal constants $c$ and $d$ such that for any solvable group $G$, any
irreducible characters $\chi,\psi$ and $\theta$ such that $\theta$ is a constituent of 
$\chi\psi$, we have that $\dl(\Ker (\theta)/(\Ker(\chi)\cap\Ker(\psi)))\leq c\eta(\chi\psi)+d$.

 In Theorem A of \cite{derived1} is proved that given any supersolvable group $G$ and any conjugacy class $A$,
 we have that $\dl(G/{\bf C}_G(A))\leq 2\eta(AA^{-1})-1$. We conjecture in \cite{derived1} that
 there exist universal constants $q$ and $r$ such that for any solvable group 
 $G$ and any conjugacy class $A$ of $\,G$, we had that 
 $\dl (G/{\bf C}_G(A))\leq q\eta(AA^{-1})+r$. In light of Theorem B and because of the 
 ``dual" situation with characters, namely  Theorem A of \cite{derivedcharacters2}, we wonder the 
 following 
 
 \begin{conj}  There exist universal constants $r$ and $s$ such that for any solvable group 
 $G$, any conjugacy classes $A$, $B$ and $D$ of $\,G$ such that $D\subseteq AB$, 
 we have that 
 $$\dl ({\bf C}_G(D)/({\bf C}_G(A)\cap{\bf C}_G(B)))\leq r\eta(AB)+s.$$
\end{conj}

We will show in Theorem \ref{equivalent} that the previous conjecture has an affirmative answer if and only 
if Conjecture of \cite{derived1} has an affirmative answer. We would like to point
out that there are several examples of ``dual results" between products of conjugacy classes
and products of character. For example,  see \cite{arad}, 
 \cite{productchar1} and \cite{productcp2}, \cite{productchar2} and \cite{productc}.
  However not every result in products of characters has a ``dual"
result in conjugacy classes, see for instance Section 3.

 In Section 4 
 we provide an example of a property in a conjugacy class  $A$ of $G$ that bounds the nilpotent class
  of  $G/{\bf C}_G(A)$, and therefore it bounds the derived length of that section.
\end{section}

\begin{section}{Proofs}
{\bf Notation.} Let $G$ be a group, $X$ be a subset of $G$ and $a\in G$. 
 Set $[a,X]=\{[a,x]\mid x\in X\}$. 
Observe that $a^X=\{a[a,x]\mid x\in X\}=a[a,X]$.

\begin{proof}[Proof of Theorem A]
Write $N={\bf C}_{G}(AB)$ and $C={\bf C}_{G}(N)$ so $AB\subseteq C$. Let $a\in A$,
$b\in B$ and $g\in C$. Then $ab$ and $ab^g$ lie in $AB$, and so lie in $C$, and since
$C$ is a subgroup, it follows that $b^{-1}b^g$ also lies in $C$. It follows then
than working in 
$G/C$, $b$ is central, so $[\left\langle b\right\rangle, G]\subseteq C$ and $[\left\langle b\right\rangle, G, N]=1$. In particular, 
$[\left\langle b\right\rangle, N,N]=1$ and so by the three-subgroups lemma, $[N',\left\langle b\right\rangle]=1$. Since this holds
for all $b\in B$, we have that $N'\subseteq {\bf C}_G(B)$. Similarly, $N'\subseteq
{\bf C}_G(A)$ and the result follows.
\end{proof}

\begin{example}\label{bunchextraspecials}
Let $m>0$ be an integer. Write $m=\prod_{i=1}^l p_i$, where $p_i$ are primes not necessarily distinct. Let $P_i$ be nonabelian of order $p_i^3$, and let 
$G$ be the direct product of the groups $P_i$. Choose noncommuting elements $a_i$ and
$b_i$ in $P_i$ and write 
$a=\prod a_i$ and 
$b=\prod b_i$. It is easy to check that $A=a^G=aZ$, $B=b^G=bZ$, where 
$Z={\bf Z}(G)$. Then $AB=abZ=D$, where $D=(ab)^G$. Also ${\bf C}_G(A)\cap {\bf C}_G(B)=Z$,
which has order $\prod p_i$, while 
${\bf C}_G(D)$ has order $\prod (p_i)^2$. Then 
$|{\bf C}_G(D)/({\bf C}_G (A)\cap {\bf C}_G(B))|=\prod_{i=1}^l p_i=m$, as wanted.
\end{example}

\begin{definition}
Let $F(n)$ be a nondecreasing function defined on the natural numbers. A group $G$ is {\bf good} for 
$F$ if for every conjugacy class $A$ of $G$, the group 
$G/{\bf C}_G(A)$ is solvable with derived length at most $F(\eta(AA^{-1}))$.
\end{definition}

In Theorem A of \cite{derived1}, it is proved that $F(n)=2n-1$ is good for all 
supersolvable groups.

\begin{theo}\label{equivalent}
Suppose that a function $F$ is good for all homomorphic images of $G$. Let $A$, $B$ and 
$D$ be conjugacy classes of $G$ with $D\subseteq AB$. Then 
${\bf C}_G(D)/ ({\bf C}_G(A)\cap {\bf C}_G(B))$ is solvable with derived length
at most $1+F(\eta(AB))$. 
\end{theo}
\begin{proof}
Let $N={\bf C}_G(D)$, and observe that ${\bf C}_G(A)\cap {\bf C}_G(B)\subseteq N$. It suffices to show that both $N/(N\cap {\bf C}_G(A))$
and $N/(N\cap {\bf C}_G(B))$ have derived length at most $1+F(\eta(AB))$. We will prove
the required length inequality for $N/(N\cap {\bf C}_G(A))$; the other inequality 
follows similarly.

Let $C={\bf C}_G(N)$ and, using the standard bar convention, write $\bar{G}=G/C$.
Since $D\subseteq C$, we see that the identity is an element of $\bar{A}\bar{B}$.
Also $\bar{A}$ and $\bar{B}$ are conjugacy classes of $\bar{G}$, so they must
be inverse classes. By Lemma 2.5 of \cite{productc} we have that
$\eta(\bar{A}\bar{B})\leq \eta(AB)$, and so  by hypothesis we have then
\begin{equation}\label{equationkey}
\dl (\bar{G}/{\bf C}_{\bar{G}}(\bar{A}))\leq F(\eta(\bar{A}\bar{B}))\leq F(\eta(AB)).
\end{equation}
Thus $\dl (G/K)=  F(\eta(\bar{A}\bar{A}^{-1}))$, where $K$ is the preimage in $G$ of 
${\bf C}_{\bar{G}}(\bar{A}))$. It follows then that 
$\dl(N/(N\cap K))\leq F(\eta(\bar{A}\bar{A}^{-1}))$. Now ${\bf C}_G(A)\subseteq K$ and so
it will be enough to show that 
$N\cap K/(N\cap {\bf C}_G(A))$ is abelian, yielding the desired derived length bound for
$N/(N\cap {\bf C}_G(A))$.

Now $K$ centralizes $A$ modulo $C$, so $[\left\langle a\right\rangle, K]\subseteq C$ and $[\left\langle a\right\rangle, K,N]=1$ for 
 $a\in A$. Then $[\left\langle a\right\rangle, (N\cap K), (N\cap K)]=1$ and by the three-subgroups lemma,
$(N\cap K)'$ centralizes $a$, and we have  $(N\cap K)'\subseteq N\cap {\bf C}_G(A)$ as
wanted.
\end{proof}

Since $F(n)=2n-1$ is good for supersolvable groups (Theorem A of \cite{derived1}), 
Theorem A follows from the previous result.
\end{section}
\begin{section}{Conjugacy class sizes}

Fix a prime $p$. Let $G$ be a finite $p$-group
 and $A$ be conjugacy classes of $\,G$.  
 In Theorem A of \cite{productc}, we proved
 that if $|A|=p^n$ for some integer $n$, then 
 $\eta(AA^{-1})\geq n(p-1)+1$. Thus, in the particular case that $G$ is $p$-group, there
 is a relation between the size of $A$ and $\eta(AA^{-1})$. 
  We want to point out a ``dual" result
 in character theory, where the ``dual"
  of a conjugacy class is an irreducible complex character $\chi$ and 
 the ``dual" of the inverse of a conjugacy class is the complex conjugate character $\overline{\chi}$
 of $\chi$, where $\overline{\chi}(g)=\overline{\chi(g)}$ for all $g\in G$. More specifically,
in Theorem A of \cite{productchar2} is proved that if $G$ is a $p$-group,
  $\chi$ is an irreducible character of 
 $G$ and $\chi(1)=p^n$, then the product $\chi\overline{\chi}$ of $\chi$ and 
 $\overline{\chi}$ has at
 least $2n(p-1)+1$ distinct irreducible constituents, i.e $\eta(\chi\overline{\chi})\geq 2n(p-1)+1$. 
 
  In Theorem A of \cite{length1} is proved that if  
 $\chi$  is an irreducible 
  character of a solvable group $G$ with $\chi(1)>1$, then $\chi(1)$ has at most $\eta(\chi\overline{\chi})-1$ different prime factors. If, in 
 addition, $G$ is supersolvable, then $\chi(1)$ has at most
 $\eta(\chi\overline{\chi})-2$ prime factors. Is there any ``dual" result in conjugacy classes
 as Theorem A of \cite{length1} in characters?  In other words, 

\begin{question}
    Does it exist a function $f:Z\rightarrow Z$ such that for 
   any solvable group $G$ and any conjugacy class $A$ of $\,G$, we have that
    $|A|$ has at most 
 $ f(\eta(AA^{-1}))$ different prime
 factors? 
\end{question}

 The answer is no, such function can not exist. More specifically, 
 
 \begin{example}
 Let $p$ be a prime and $P$ be a group of order $p$. 
 Let $G$ be the group of order $p(p-1)$, where $P\triangleleft G$ 
 and $G/P$ induces on $P$ the full group of automorphisms of $P$.
 Then $P$ contains just one nontrivial conjugacy class $A$ of $G$, namely $A=P\setminus \{1\}$. Observe that $AA^{-1}=P$ and 
 so $\eta(AA^{-1})=2$. Also $P={\bf C}_G(A)$, and thus 
 $|G/{\bf C}_G(A)|=p-1$. This is obviously unboundedly large, and by a
 result of Erdos \cite{erdos}, it has unboundedly many prime factors.
 \end{example}  
 
{\bf Remark.} Let $G$ be the group as in the previous Example. Let $\lambda$ be
an irreducible character of $N$. Then $\lambda$ is a linear character and the induced
character $\lambda^G$ is an irreducible character of degree $p-1$. Set $\chi=\lambda^G$.
We can check then  $\eta(\chi\overline{\chi})=p$, namely the irreducible constituents
of $\chi\overline{\chi}$ are $\chi$ and all irreducible character with kernel containing
$N$. Since $n-1\leq 2^{n-2}$ for any integer $n>0$, then $p-1$ has at most $p-2$ distinct prime factors. 
\end{section}

\begin{section} {Conjugacy classes and nilpotent class}
 Let ${\bf Z}_1(G)={\bf Z}(G)$ be the center of the group $G$ and by induction define
  the $i$-center of $G$ as
  ${\bf Z}_i(G)/{\bf Z}_{i-1}(G)= {\bf Z}(G/{\bf Z}_{i-1}(G))$.
  The following is a well known result.
  
  \begin{lem}
Let $N$ be a group. Write $Z_m=Z_m(N)$ for the $m$-th center of $N$ and write $N^m$ for the 
$m$-th term of the lower central series of $N$. Then $[N^m, Z_m]=1$. 
\end{lem}
\begin{proof}
Induct on $m$. For $m=1$, we have $[N^1, Z_1]=[N,{\bf Z}(N)]=1$, as needed. For $m>1$, we want $[N^{m-1},N,Z_m]=1$.
We have  $[N, Z_m, N^{m-1}]\subseteq [Z_{m-1}, N^{m-1}]=1$ by the inductive hypothesis. Now work in $\bar{N}=N/Z$ where $Z=Z_1={\bf Z}(N)$. Note that 
$\overline{N^{m-1}}=(\overline{N})^{m-1}$ and $\overline{Z_m}=Z_{m-1}(\bar{N})$. Then 
\begin{equation*} 
1=[Z_{m-1}(\bar{N}),(\bar{N})^{m-1}]=[\overline{Z_m}, \overline{N^{m-1}}]=[\overline{Z_m, N^{m-1}}],
\end{equation*}
\noindent and we have $[Z_m, N^{m-1}]\subseteq Z$. Then $[Z_m, N^{m-1},N]=1$, and the result follows by the
three-subgroups lemma. 
\end{proof}

\begin{theo}
Let $N\trianglelefteq G$ and $a \in G$. Assume that $[N,a]\subseteq Z_m$, where $Z_m=Z_m(G)$ is the $m$-th center of $N$. Let
$C={\bf C}_N(a^N)$. Then $N/C$ is nilpotent of class at most $m$.
\end{theo}
\begin{proof}
Induct on $m$. if $m=0$ we are assuming that $[N,a]=1$ and so $N$ centralizes all of $a^N$ and $C=N$. In this case
$N/C$ is nilpotent of class zero. We can assume therefore that $m>0$. Our goal is to show that
$N^{m+1}\subseteq C$.

Now let $x\in a^N$. We want $[N^{m+1}, \left\langle x\right\rangle ]=1$, or equivalently, $[N^m,N,\left\langle x\right\rangle]=1$. Now $[N,\left\langle x\right\rangle, N^m]\subseteq[Z_m, N^m]=1$ by the previous lemma.

Let $Z=Z_1={\bf Z}(N)$ and write $\bar{G}=G/Z$. Then 
\begin{equation*}
Z_{m-1}(\bar{N})=\overline{Z_m} \supseteq\overline{[N,a]}=[\bar{N}, \bar{a}],
\end{equation*}
\noindent and thus by the inductive hypothesis, $N/B$ is nilpotent of class at most $m-1$, where $B$ is the preimage
in $N$ of ${\bf C}_{\bar{N}}((\bar{a})^{\bar{N}})={\bf C}_{\bar{N}}(\overline{a^N})$. In particular,
$N^m\subseteq B$ so $x$ centralizes $N^m$ modulo $Z$, and we have $[\left\langle x\right\rangle, N^m]\subseteq Z$, and
hence $[\left\langle x\right\rangle, N^m, N]=1$. The three-subgroups lemma now yields $[N^m, N, \left\langle x\right\rangle]=1$, as wanted.
\end{proof}
{\bf Remark.} Let $G_1=C_2 \wr C_2$ be the wreath product of $C_2$  by $C_2$, where 
$C_2$ is the cyclic group of order 2. Thus $|G_1|=8$ and $G_1$ is non abelian.
Let $a_1\in G_1\setminus {\bf Z}(G_1)$. We can check that ${a_1}^{G_1}=a{\bf Z}(G_1)$
and $G_1/{\bf C}_{G_1}(a_1^{G_1})$ is abelian and so it is nilpotent of class 1.
    
    Let $N=G_1\times G_1$ and $a_2=(a_1, a_1)$ in $N$. Observe that $C_2$ acts on $N$ 
    by permuting the entries. Set $G_2=C_2 N$. We can check that 
    $a_2^{G_2}\subset a_2 {\bf Z}_2 (G_2)$ but $a_2^{G_2}\not\subset a_2 {\bf Z} (G_2)$,
    and  $G_2/({\bf C}_{G_2}(a_2^{G_2}))$ is nilpotent of class 2.
    
    The author wonders if given any integer $m>2$, we can find a group $G$ with an
    element $a\in G$ such that $a^G \subseteq a {\bf Z}_m (G)$ 
    and $G/{\bf C}_{G}(a^{G})$ is nilpotent of class $m$.
    
      {\bf Acknowledgment.} I thank Harvey Blau and Michael Bush for useful discussions and emails. 
   I would like to thank  the Instituto de Matematicas, UNAM, unidad Morelia, and 
   Gerardo Raggi for their hospitality while visiting. Also, I thank the referee of a previous version of this note
   for helpful comments and suggestions of how to improve both the presentation and the results 
   of this note;
   those are included here with her/his permission. 
\end{section}


\begin{thebibliography}{9}
\bibitem{length1} E. Adan-Bante, Products of characters and derived length, J. of Algebra, 266 (2003), 305-319.

\bibitem{productchar1} E. Adan-Bante, Products of characters and finite $p$-group, J. of Algebra, 277 (2004), 236-255.

\bibitem{productchar2} E. Adan-Bante,
 Products of characters and finite p-groups II, Archiv der Mathematik, 82  No 4 (2004), 289-297.

\bibitem{derivedcharacters2} E. Adan-Bante, Products of characters and derived length II, J. Group Theory, 
8 (2005), 453-459.

\bibitem{productc} E. Adan-Bante, Conjugacy classes and finite p-groups, Arch. Math. 85 (2005) 297-303.

\bibitem{producth} E. Adan-Bante, Homogeneous products of conjugacy classes, Arch. Math. 86 (2006)
289-294.

\bibitem{derived1} E. Adan-Bante, Derived length and products of conjugacy classes, 
Israel J. Math. 168 (2008), 93-100.

 
 

\bibitem{productcp2} E. Adan-Bante, On nilpotent groups and conjugacy classes, to appear 
Houston Journal of Mathematics. 


\bibitem{arad}
Z. Arad, E. Fisman,
An analogy between products of two conjugacy classes and products 
of two irreducible characters in finite groups, Proc. of the  Edinburgh Math. Soc. 30 (1987), 7-22. 
 
 \bibitem{erdos}P.  Erdos, On the normal number of prime factors of p-1 and some related problems concerning Euler's $\phi$-function,
Q. J. Math., Oxf. Ser. 6, 205-213 (1935). 

\end{thebibliography}
\end{document}